\documentclass[a4paper,11pt]{article}
\usepackage{amssymb}
\usepackage{amsthm}
\usepackage{amsmath}
\usepackage{latexsym}
\usepackage{epsfig}
\usepackage{amscd}
\usepackage{graphicx}
\usepackage[dvips]{color}
\newtheorem{theorem}{Theorem}[section]
\newtheorem{lem}[theorem]{Lemma}
\newtheorem{cor}[theorem]{Corollary}

\theoremstyle{definition}
\newtheorem{definition}[theorem]{Definition}

\theoremstyle{remark}
\newtheorem{rem}[theorem]{Remark}

\numberwithin{equation}{section}

\title{Rigidity of Conformally Compact Manifolds with the Round Sphere as the Conformal Infinity}
\author{ Satyaki Dutta }

\begin{document}
\maketitle

\begin{abstract}
In this paper we prove that under a lower bound on the Ricci curvature and an asymptotic assumption on the scalar curvature, a complete conformally compact manifold $(M^{n+1},g)$, with a pole $p$ and with the conformal infinity in the conformal class of the round sphere, has to be the hyperbolic space. 
\end{abstract}

\section{Introduction}

Rigidity problems for conformally compact manifolds$(M^{n+1}, g)$ have been studied by several mathematicians. Several partial results have been proved over the years. It has been shown in \cite{min} that a strongly asymptotic hyperbolic spin manifold of dimension $n$ and with scalar curvature $R \geq -n(n+1)$ has to be the hyperbolic space. Leung proved in \cite{leu} that any conformally compact Einstein manifold, that is asymptotically hyperbolic of order greater than 2 is in fact hyperbolic. In \cite{qin}, J.Qing showed  that a conformally compact Einstein manifold of dimension $ n \leq 7$, and with round sphere as it's conformal infinity must be the hyperbolic space. His proof relies on the positive mass theorem of Schoen and Yau, hence the constraint on the dimension. Anderson showed in \cite{ma4} that a conformally compact Einstein manifold with conformal infinity the round sphere is the hyperbolic space he assumes no restriction on the dimension. We prove here a similar result, where the manifold is "close to" but not Einstein and with no restriction on the dimension.  Recently Shi and Tian in \cite{shitian} proved that an ALH manifold with a pole and with sectional curvatures going to -1 at a rate, $|K_{ij}+1| = O(e^{-\alpha r})$, where $\alpha > 2$ and $Ric \geq -ng$ has to be the hyperbolic space. We prove the same conclusion, but with a weaker assumption on the sectional curvatures, that is $|K_{ij}+1| = O(\phi(r)e^{-2r})$, where $\phi(r) \geq 0$ and $\phi \in L^{1}(\mathbb{R}^{+})$ . 

Before stating our main result, we briefly review the definition of $C^{m,\alpha}$ conformally compact manifolds. 

\begin{definition}
A complete manifold $(M^{n+1},g)$ is said to be $C^{m,\alpha}$ conformally compact if the conformally equivalent metric $\bar{g} = \rho^{2}g$ extends to a $C^{m,\alpha}$ metric on $\partial M$. Where $\rho$ is a $\emph{defining function}$, i.e, $\rho : \bar{M} \rightarrow \mathbb{R}$, such that $\rho > 0$ in $M$, $\rho|_{\partial{M}} = 0$, and $d\rho \neq 0$ on $\partial M$. 
\end{definition}
Additionally, if $|\bar{\nabla} \rho|_{\bar{g}} = 1$ in a neighborhood of $\partial M$, we say that $\rho$ is a $\emph{geodesic defining function}$. A simple and interesting fact, shown in \cite{ma2}, is that if the compactification is atleast $C^{2}$, then for a given boundary metric, every such conformal compactification has a unique geodesic defining function. Defining functions are unique up to multiplication by positive  functions on $\bar{M}$. Therefore, only the conformal class $[\bar{g}]$ is uniquely determined by $g$, same goes with the conformal class $[\gamma]$ of the induced boundary metric. The $C^{m,\alpha}$ $\textbf{Conformal Infinity}$ is defined as the $C^{m,\alpha}$ conformal class of metrics on the boundary.
 Its easy to see that for conformally compact manifolds 
\begin{equation}
\label{curv}
 |K_{ij} + 1| = O(\rho^{2}). 
\end{equation} 
 That is, such manifolds are asymptotically hyperbolic. 
We now state our main result in the following Theorem:
\begin{theorem}

Let $(M^{n+1}, g)$, $n \geq 3$ be a complete conformally compact manifold with a pole $p$ with conformal infinity $(S^{n}, [\gamma_{o}])$ and $Ric \geq -ng$, where $[\gamma_{0}]$ is the conformal class of the round metric. Let $t(x)= dist_{g}(p,x)$ and $|R +n(n+1)| = o(e^{-2t})$, where $R$ is the scalar curvature of $(M^{n+1}, g)$. Then $(M^{n+1}, g)$ is $(B^{n+1}, g_{-1})$. 

\end{theorem}
In contrast to earlier results, the Theorem requires only a lower bound on the Ricci curvature. We make no additional assumption on the sectional curvatures other than what is inherited through conformally compactness. The scalar curvature of any conformally compact manifold satisfies $|R +n(n+1)| = O(e^{-2t})$, by virtue of (\ref{curv}). We require a slightly stronger decay rate. This also relates to the gap Theorems of \cite{gw}, \cite{ks}, that assert that there are gaps in the positive and negative sides of space forms in the set of Riemannian manifolds. Throughout the article we stick with the convention that the scalar curvature is the trace of the Ricci curvature.

The organization of the paper is as follows. We start off with a geodesic defining function $\rho: \bar{M} \to \mathbb{R}$ of the boundary that is $|\bar{\nabla} \rho|_{\bar{g}} \equiv 1$ in a neighborhood of the boundary of $M$. Set $\rho = 2e^{-r}$, it is easy to show that $|\nabla r|_{g} = 1$ if and only if $|\bar{\nabla} \rho|_{\bar{g}} = 1$. Therefore, outside of a compact set $K \subset M$, $|\nabla r|_{g} = 1$, and therefore $r$ is a distance function. Choose $r_{0}$, such that $M_{r_{0}} = \{ x \in M | \rho(x) > \rho_{0} = 2e^{-r_{0}} \} $  contains the cut locus of $\rho$, the set $K$ and the pole $p$. For any $x$ in $M$ and outside of $M_{r_{0}}$, define $t(x)= dist_{g}(p,x)$. Recall, $p$ is the pole and therefore $t$ is a smooth distance function on $M$. Clearly the function $c(x) = t(x) - r(x)$ is a smooth bounded function in $M - M_{r_{0}}$.
Next we look at two compactifications. $\bar{g} = 4e^{-2r}g$ and $\bar{g}' = 4e^{-2t}g$. By our assumption, $(\bar{M}, \bar{g})$ is the smooth compactification of $(M,g)$, with the round sphere as its boundary. We restrict the function $c: M \rightarrow \mathrm{R} $ to the level-sets of $t$, $c_{t}: (\bar{\Sigma}'_{t}, \bar{g}'|_{\bar{\Sigma}'_{t}}) = (\bar{\Sigma}'_{t}, \bar{\gamma}'_{t}) \rightarrow \mathrm{R}$.  We then find an uniform bound on $||c_{t}||_{C^{1}} $ for all $t$ sufficiently large. Using Arzela-Ascoli's Theorem we conclude that $c$ has a $C^{0,1}$ extension to the boundary. And hence $\bar{\gamma}'_{t}$ extends to a metric $\gamma$ that is $C^{0,1}$ conformal to $\gamma_{0}$. Next we show that the scalar curvatures $\bar{R}'_{t}$, on the level sets, $(\bar{\Sigma}'_{t}, \bar{\gamma}'_{t})$, are less than $n(n-1) + o(1)$, and the volumes of the level sets are bounded above by the volume of the unit round sphere. We use these estimates and the Yamabe quotient for $\gamma$, and show that $\gamma$ infact has to be smooth. Thereafter we apply Obata's Theorem, and conclude that $(\partial M, \gamma)$ is isometric to $(S^{n}, \gamma_{0})$. Finally, we complete the proof of the Theorem by applying Bishop-Gromov Volume comparison Theorem on balls centered at the pole $p$.

This work formed the author's PhD dissertation. It is a pleasure to thank Michael Anderson for introducing me to this topic and for his guidance and inspiration all along. I would also like to thank Marcus Khuri for useful conversations on several occasions.
 
\section{Asymptotic Bounds on the Shape Operators}

It is well known that the shape operator $S$, satisfies the following,
\begin{equation}\label{shop1}
\partial_{r}(S_{i}^{j}) +(S_{i}^{k})(S_{k}^{j}) = -R_{i}^{j}
\end{equation}
Taking its trace we have,
\begin{equation} \label{shop2}
\partial_{r}m + tr\Big(S_{i}^{k}S_{k}^{j} \Big) = -Ric(\partial_{r}, \partial_{r}) 
\end{equation}
where $m = tr(S_{i}^{j})$. We will use these equations several times in this article. We will now find asymptotic bounds on the shape operators $\nabla^{2} r$ and $\nabla^{2} t$ with respect to the complete metric $g$. 

\begin{lem}\label{lem1}

 Let $(M,g)$ be a complete conformally compact manifold. Let $r$ and $t$ be smooth distance functions as defined earlier. For all $r$ and $t$ sufficiently large, the eigenvalues $\{\lambda_{i}(r) \}$ and $ \{ \mu_{i}(r) \} $ of the shape operators $\nabla^{2}r$ and $\nabla^{2}t$, respectively are non-negative and bounded from above.

\end{lem}

\begin{proof}
Since $(M,g)$ is conformally compact, we have $|K_{ij} + 1| = O(e^{-2r})$.  So outside of a compact set $K$, we can assume that $ -\alpha^{2} \leq K_{ij} \leq -\frac{1}{\alpha^{2}} $, for $ \alpha > 1$ but close to 1. Without loss of generality, we can assume $K$ to be the same compact set as above. 
Choose $r'_{0}$ such that $K$ lies inside $M_{r'_{0}}$. For any $r > r'_{0}$ we have the following differential inequality for the each eigenvalue, $\lambda_{i}(r)$ of the shape operator $\nabla^{2}r$ as in ~(\ref{shop1}),

$$ \frac{1}{\alpha^{2}} \leq \lambda'_{i}(r) + \lambda^{2}_{i}(r) \leq \alpha^{2}. $$
From the second inequality, we get,
$$ \frac{\lambda'_{i}(r)}{\alpha^{2} - \lambda^{2}_{i}(r)} \leq 1. $$
Integrating from $r'_{0}$ to $r$, we get,

$$ \ln \Big (\frac{\alpha + \lambda_{i}(r)}{\alpha - \lambda_{i}(r)}\frac{1}{\hat{C}_{0}(r'_{0})}\Big ) \leq 2 \alpha (r - r'_{0}), $$
where $\hat{C}_{0}(r'_{0}) = \ln(\frac{\alpha + \lambda_{i}(r'{0})}{\alpha - \lambda_{i}(r'_{0})}) $
\newline 
or,

 $$ \lambda_{i}(r) \leq \alpha \Big (\frac{C_{0}(r'_{0})e^{2 \alpha r} - 1}{C_{0}(r'_{0})e^{2 \alpha r} + 1}\Big ), $$
where $C_{0}(r'_{0}) = \hat{C}_{0}(r'_{0})e^{-2\alpha r'_{0}} $.
The first inequality, 

 $$ \frac{1}{\alpha^{2}} \leq \lambda'_{i}(r) + \lambda^{2}_{i}(r) $$
gives,

$$ \lambda_{i}(r) \geq \frac{1}{\alpha} \Big (\frac{C'_{0}(r'_{0})e^{\frac{2}{\alpha} r} - 1}{C'_{0}(r'_{0})e^{\frac{2}{\alpha} r} + 1}\Big ), $$
where $C'_{0}(r'_{0}) = \ln(\frac{\frac{1}{\alpha} + \lambda_{i}(r'{0})}{\frac{1}{\alpha} - \lambda_{i}(r'_{0})})e^{-\frac{2}{\alpha} r'_{0}} $
Without loss of generality, we can assume $r_{0}' = r_{0}$. 
Therefore, we see that $ 0 \leq \lambda_{i}(r) \leq C_{1}(r_{0}) $  for all $ r > r_{0}$. 
Since for any point $x \in M$, $|t(x) - r(x)|$ is bounded, we also have, $|K_{ij} + 1| = O(e^{-2t})$.
Therefore, a similar analysis for the eigenvalues $\mu_{i}(t)$ of $\nabla^{2}t$ gives us the following inequality 

$$  0 \leq \mu_{i}(t) \leq E_{1}(t_{0}) $$
for every $\mu_{i}(t)$ and for all $t > t_{0}$. 

\end{proof}
The next Lemma shows that the shape operators actually tend to 1 at a certain rate. Under a stronger assumption on the Ricci curvature, we will show that the Laplacian of $t$ satisfies certain asymptotic bound. We wish to emphasize here that this extra Ricci curvature assumption is not necessary to prove our result and is not used anywhere in the proof.

\begin{lem} \label{lem2}
Let $(M,g)$ be a complete conformally compact manifold. Let $r$ and $t$ be smooth distance functions as defined earlier. For $r > r_{0} $ and $t > t_{0}$, we have:

\begin{eqnarray}
\label{eq3.1.1}
\nabla^{2}r(e_{i},e_{j}) = \delta_{i}^{j} + o(e^{-\frac{3}{2}\beta r}) \\
\label{eq3.1.2}
\nabla^{2}t(e'_{i},e'_{j}) = \delta_{i}^{j} + o(e^{-\frac{3}{2}\beta t}) 
\end{eqnarray}
where $\beta < 1 $,  $ \{ e_{i} \}_{i=1}^{n} $ are tangent to the level sets $\Sigma_{r}$ and $ \{ e'_{i} \}_{i=1}^{n} $ are tangent to the level sets $\Sigma_{t}$.
If, in addition $ -(n-\psi(t)e^{-2t})g \leq Ric \leq -ng $, where $\psi(t) \geq 0 $ and $\psi(t) \in L^{1}(\mathbf{R}^{+})$, then 
\begin{equation}
\label{eq3.1.3}
\Delta t = n + O(e^{-2t})
\end{equation} 

\end{lem}

\begin{proof}

For $r > r_{0}$, we have the following inequality for $\{ \lambda_{i}(r) \}_{i=1}^{n}$ from the proof of Lemma ~\ref{lem1}:

$$ 0 \leq \lambda_{i}(r) \leq C_{1}(r_{0}). $$
We also have:

  $$ -1 - Ce^{-2r} \leq K_{ij} \leq -1 + Ce^{-2r}. $$
This gives us the following differential inequality:

\begin{equation}
\label{eq3.1.4}
1 - Ce^{-2r} \leq \lambda'_{i}(r) + \lambda^{2}_{i}(r) \leq 1 + Ce^{-2r}.
\end{equation}
Let $u_{i}(r) = \lambda_{i}(r) - 1 $, therefore $ -1 \leq u_{i}(r) \leq C_{1}(r_{0}) $.
Thus, 
\begin{eqnarray}
(u^{2}_{i}(r))' + 2u^{2}_{i}(r) & \leq & (u^{2}_{i}(r))' + (4 + 2u_{i}(r))u^{2}_{i}(r)  \nonumber \\
 & \leq & 2Cu_{i}(r)e^{-2r}. \nonumber
\end{eqnarray} 
That is,

 $$ (u^{2}_{i}(r))' + 2u^{2}_{i}(r) \leq 2Cu_{i}(r)e^{-2r}. $$
Choose $0 < \beta' < 1 $, and write the above inequality in the following manner,

\begin{eqnarray}
(u^{2}_{i}(r))' + 2\beta'u^{2}_{i}(r) & \leq & (u^{2}_{i}(r))' + 2u^{2}_{i}(r) \nonumber \\ 
 & \leq & 2Cu_{i}(r)e^{-2\beta'r}e^{-2(1-\beta')r}. \nonumber
\end{eqnarray}
Therefore, bringing $e^{-2\beta'r}$ over to the other side,

  $$ (u^{2}_{i}(r)e^{2\beta'r})' \leq 2Cu_{i}(r)e^{-2(1-\beta')r}. $$
Integrating from $r_{0}$ to $r$,

$$ u^{2}_{i}(r)e^{2\beta'r} \leq C_{2}^{2}(r_{0}), $$
 or,

\begin{equation}
\label{eq3.1.5}
 |u_{i}(r)| \leq C_{2}(r_{0})e^{-\beta'r}.
\end{equation}
This estimate is far weaker than what we had set out to show. To improve this we use the following differential inequality,

\begin{eqnarray}
\label{eq3.1.6}
(u^{2}_{i}(r))' + (4 - 2|u_{i}(r)|)u^{2}_{i}(r) & \leq & (u^{2}_{i}(r))' + (4 + 2u_{i}(r))u^{2}_{i}(r)  \nonumber \\
 & \leq & 2C|u_{i}(r)|e^{-2r}.
\end{eqnarray}
Putting back the estimate on $|u_{i}(r)|$ obtained from (~\ref{eq3.1.5}) to (~\ref{eq3.1.6})

 $$ (u^{2}_{i}(r))' + 4u^{2}_{i}(r) \leq 2C.C_{2}(r_{0})e^{-(2+\beta')r} + 2C_{2}^{3}(r_{0})e^{-3\beta'r}. $$
Choose $ 0 < \beta < \beta' < \frac{(2+\beta')}{3} $, and split the right hand side of the inequality in the following way,

$$ (u^{2}_{i}(r))' + 4u^{2}_{i}(r) \leq e^{-3\beta r}[2C.C_{2}(r_{0})e^{-(2+\beta'-3\beta)r} + 2C_{2}^{3}(r_{0})e^{-3(\beta' - \beta)r}] $$
or,

\begin{eqnarray}
(u^{2}_{i}(r))' + 3\beta u^{2}_{i}(r) & \leq & (u^{2}_{i}(r))' + 4u^{2}_{i}(r)  \nonumber \\
 & \leq & e^{-3\beta r}[2C.C_{2}(r_{0})e^{-(2+\beta'-3\beta)r} + 2C_{2}^{3}(r_{0})e^{-3(\beta' - \beta)r}] \nonumber
 \end{eqnarray}
or,

 $$(u^{2}_{i}(r)e^{3\beta r})' \leq [2C.C_{2}(r_{0})e^{-(2+\beta'-3\beta)r} + 2C_{2}^{3}(r_{0})e^{-3(\beta' - \beta)r}] $$
Integrating from $r_{0}$ to $r$,

 $$ u^{2}_{i}(r)e^{3\beta r} \leq C_{3}^{2}(r_{0}) $$
or,
\begin{equation}
\label{eq3.1.7}
 \lambda_{i}(r) = 1 + o(e^{-\frac{3}{2}\beta r}). 
\end{equation}
Similar analysis with the eigenvalues $\mu_{i}(t)$ proves the other claim.
\newline
Here we prove our last claim. As before, let $\mu_{i}(t)$ 's be the eigenvalues of the Hessian of t, therefore, $\Delta t = \Sigma_{i=1}^{n} \mu_{i}(t)$.  
Let us assume $ \mu_{i}(t) = 1 + T_{i}(t)e^{-2t} $, we would show that $|\Sigma_{i=1}^{n}T_{i}(t)| = O(1) $. From ~(\ref{eq3.1.2}) we already know that for each $i$, $|T_{i}(t)| \leq O(e^{(2-\frac{3}{2}\beta) t}).$ Furthermore, 
by ~(\ref{shop2}), we get,
 
  $$ n - \psi(t)e^{-2t} \leq (\Sigma_{i=1}^{n} \mu_{i}(t))' + \Sigma_{i=1}^{n} \mu^{2}_{i}(t) \leq n. $$ 
Replacing $\mu_{i}(t)$ by $1+T_{i}(t)e^{-2t}$, we get,
$$ -\psi(t) \leq (\Sigma_{i=1}^{n}T_{i})' + (\Sigma_{i=1}^{n}T_{i}^{2})e^{-2t} \leq 0. $$
Now, the second inequality, that is  $(\Sigma_{i=1}^{n}T_{i})' + (\Sigma_{i=1}^{n}T_{i}^{2})e^{-2t} \leq 0 $ implies, $ (\Sigma_{i=1}^{n}T_{i})' \leq 0$, that is $\Sigma_{i=1}^{n}T_{i}(t)$ is non-increasing.
Therefore it can only go to $ - \infty $. Whereas the first inequality stops that from happening.
\begin{eqnarray}
  -\psi(t) & \leq & (\Sigma_{i=1}^{n}T_{i})' + (\Sigma_{i=1}^{n}T_{i}^{2})e^{-2t}, \\  \nonumber
  -\psi(t) - (\Sigma_{i=1}^{n}T_{i}^{2})e^{-2t} & \leq & (\Sigma_{i=1}^{n}T_{i})'. \nonumber
\end{eqnarray}
Since we already know that, for each $i$,   $|T_{i}(t)| \leq O(e^{(2-\frac{3}{2}\beta) t})$, and $2-\frac{3}{2} \beta \leq 1- \epsilon $, for some $\epsilon > 0$. We have the left hand side in $L^{1}(\mathbf{R}^{+}) $, and therefore $\Sigma_{i=1}^{n}T_{i}(t)$ is greater than some finite number. 
Hence we conclude that $|\Sigma_{i=1}^{n}T_{i}(t)| = O(1) $.

\end{proof}
\begin{rem}\label{rem3.1.1}
Since the compactification by the function $r$ is smooth, all curvature quantities, in particular $\bar{Ric}$, are bounded. A stronger asymptotic condition holds for $ \nabla^{2}r$, i.e. $\nabla^{2}r(e_{i},e_{j}) = \delta_{i}^{j} + O(e^{-2r})$. 
\end{rem}
The following corollary is a direct consequence of Lemma ~\ref{lem2}.

\begin{cor} \label{cor1}
The second fundamental form of the level-sets $( \bar{\Sigma_{t}}', \bar{\gamma}'_{t} )$,  $\bar{\nabla}'^{2}t(\bar{e}_{i}', \bar{e}_{j}') = \bar{S}'(\bar{e}_{i}',\bar{e}_{j}') \rightarrow 0$  as $t \rightarrow \infty$, where $\{ \bar{e}'_{i} \}_{i=1}^{n}$ are unit tangent vectors to $( \bar{\Sigma_{t}}', \bar{\gamma}'_{t}).$
\end{cor}

\begin{proof}
 
 Since we are considering the metric $\bar{g}' = e^{-2t}g$, where $t(x) = dist_{g}(x,p)$, the unit normal $\bar{N}'$ to the 
 level set$( \bar{\Sigma_{t}}', \bar{\gamma}'_{t} )$, in $\bar{g}'$ metric can be written as,
 
 \begin{equation}
  \bar{N}' = \bar{\nabla}'e^{-t} = -e^{-t} \bar{\nabla}' t.
 \end{equation}
 
 %$N = \nabla t$ is the unit normal to the level sets $(\Sigma_{t}, \gamma_{t})$
Therefore, the second fundamental form of the level set $( \bar{\Sigma}'_{t}, \bar{\gamma}'_{t} )$ is as follows,
 
 \begin{eqnarray}
 \bar{S}'(\bar{e}_{i}',\bar{e}_{j}') & = & \bar{g}'(\bar{\nabla}'_{\bar{e}_{i}'} \bar{N}', \bar{e}_{j}')  \nonumber\\
 & = & g(\bar{\nabla}'_{e_{i}'}\bar{N}', e_{j}')   \nonumber\\
 & = & g(\nabla_{e_{i}'}\bar{N}', e_{j}') - g(\nabla t, e_{i}')g(\bar{N}',e_{j}') - g(\nabla t,e_{j}')g(\bar{N}',e_{i}') \nonumber \\
 &   & + g(\nabla t, \bar{N}')g(e_{i}',e_{j}'). \nonumber
 \end{eqnarray}
 Since, $\{ e_{i}' \}_{i=1}^{n}$ are tangent to the level-sets and $\nabla t$ is normal, the second and the third terms above 
 are zero. Hence we have,
 
 \begin{eqnarray}
 \bar{S}'(\bar{e}_{i}',\bar{e}_{j}')  & = &  g(\nabla_{e_{i}'}\bar{N}', e_{j}') + g(\nabla t, \bar{N}')\delta_{ij} \nonumber\\
  & = & -g(\nabla_{e_{i}'} (e^{-t}\bar{\nabla}'t),e_{j}') + e^{-t}\delta_{ij}g(\nabla t, \bar{\nabla}' t) \nonumber\\
  & = & e^{-t}(g(\bar{\nabla}'t, \nabla_{e_{i}'}e_{j}')) - e_{i}'(g(e^{-t}\bar{\nabla}'t,e_{j}')) + e^{t}\delta_{ij}g(\nabla t, \nabla t). \nonumber
 \end{eqnarray}
 The second term, $ g(\bar{\nabla}'t, e_{j}') = 0$, since $\bar{\nabla}' t$ is normal to the level-sets.
 To simplify the first term we use the following,
 
 \begin{eqnarray}
 g(\bar{\nabla}'t,X) & = &  e^{2t}\bar{g}'(\bar{\nabla}'t, X) \nonumber\\
  & = & e^{2t}g(\nabla t, X). 
 \end{eqnarray} 
 Thus we have,
 
 \begin{eqnarray}
  e^{t}(g(\nabla t, \nabla_{e_{i}'}e_{j}')) + e^{t}\delta_{ij}g(\nabla t, \nabla t) & = & e^{t}(e_{i}'(g(\nabla t, e_{j}')) - g(\nabla_{e_{i}'} \nabla t, e_{j}') + \delta_{ij}) \nonumber\\
   & = & -e^{t}(g(\nabla_{e_{i}'}N, e_{j}') - \delta_{ij})  \nonumber\\
   & = & -e^{t}(\nabla^{2}t(e_{i}',e_{j}') - \delta_{ij}). \nonumber
 \end{eqnarray}
Hence, replacing $\nabla^{2}t(e_{i}',e_{j}')$ by $ \delta_{ij} + O(e^{-\frac{3}{2}\beta t}) $ proved in Lemma~(\ref{lem2}), we get,
 
 \begin{equation}
 \label{eq3.1.11}
  \bar{S}'(\bar{e}_{i}', \bar{e}_{j}') = O( (e^{-(\frac{3}{2}\beta-1)t}).
 \end{equation} 
 
 \end{proof}

 We now prove the following Lemma, which is one of the key estimates for the result in \cite{shitian}. We show here in this Lemma that the main Theorem in \cite{shitian} stays valid even under a weaker assumption on the asymptotic behavior of the sectional curvatures. One readily finds that assuming the following Lemma is true, the Theorem of \cite{shitian} follows immediately.

\begin{lem} \label{lem3}
Let $(M^{n+1},g)$, $n \geq 2$ be a complete manifold with a pole, $p$. If the sectional curvatures, $K_{ij}$ satisfy
\begin{equation}
 -1 - \phi(t)e^{-2t} \leq \mathrm{K}_{ij} \leq -1 + \phi(t)e^{-2t}
\end{equation}
where, $t(x) = dist_{g}(x,p)$ $ \phi \in \mathrm{L}^{1}(\mathbb{R^{+}})$ and $\phi \geq 0. $ and $Ric(g) \geq -ng$,
then, for $t > t_{0}$, we have:
\begin{equation}
\label{eq3.1.13}
  \nabla^{2}t(e_{i}',e_{j}')=g(S(e_{i}'),e_{j}') = \delta_{ij} + O(e^{-2t}).
\end{equation}

\end{lem}

\begin{proof}

 Our goal here is to improve the estimates derived in Lemma ~(\ref{lem2}), by using the extra assumption on the sectional curvatures. We begin with ~(\ref{eq3.1.4}) and under the stronger assumption on the sectional curvatures, we have:
 
 \begin{equation}\label{eq3.1.14}
 1-\phi(t)e^{-2t} \leq \mu_{i}'(t)+\mu_{i}^{2}(t) \leq 1+\phi(t)e^{-2t}. 
 \end{equation}
Let, 
 
  $$ v_{i}(t) = (\mu_{i}(t)-1)e^{2t}. $$
In terms of $v_{i}(t)$, ~(\ref{eq3.1.14}) becomes:
 
  $$ v_{i}^{2}(t) \leq (\phi(t) - v_{i}'(t))e^{2t}. $$
Which implies that $ \phi(t) \geq v_{i}'(t) $. Integrating this from $t_{0}$ to $t$ and using the fact that $\phi \in \mathrm{L}^{1}(\mathbb{R^{+}}) $, we get $ v_{i}(t) \leq C_{1}$, and hence 

\begin{equation}
\mu_{i}(t) \leq 1 + C_{2}e^{-2t}
\end{equation}
To get a lower bound, we write 

 $$ \mu_{i}(t) = 1 + T_{i}(t)e^{-2t} $$
We already know from ~(\ref{eq3.1.14}) and the last part of Lemma~(\ref{lem2}) with extra Ricci curvature assumption that obviously holds here, that for each $i$ and $t \geq t_{0}$,
 
 $$ -C_{3}e^{(2-\frac{3}{2} \beta)t} \leq T_{i}(t) \leq C_{2}.  \qquad \qquad  (*)$$
By ~(\ref{shop2}), we get:

  $$ (\Sigma \mu_{i}(t))' + \Sigma \mu^{2}_{i}(t) = - Ric(N,N).$$
Since $|K_{ij} + 1| \leq \phi(t)e^{-2t}$, we have a bound on $Ric(N,N)$
    $$ |Ric(N,N) + n| \leq n\phi(t) e^{-2t}. $$  
Therefore,

  $$ n - n\phi e^{-2t} \leq (\Sigma \mu_{i})' + \Sigma \mu^{2}_{i} \leq n + n\phi e^{-2t}.$$
Replacing $\mu_{i}(t)$'s by $1+T_{i}(t)e^{-2t}$'s, we get:
$$ -n\phi(t) e^{-2t} \leq (\Sigma T_{i}(t))'e^{-2t} + (\Sigma T^{2}_{i}(t))e^{-4t} \leq n\phi(t) e^{-2t}.$$
or,

\begin{equation}
  -n\phi(t) \leq (\Sigma T_{i}(t))' + (\Sigma T^{2}_{i}(t))e^{-2t} \leq n\phi(t). 
\end{equation}
Now, since $T^{2}_{i}(t) \leq max\{C^{2}_{2}, C^{2}_{3}e^{(4 - 3\beta)t} \} $. We have 
$$T^{2}_{i}(t)e^{-2t} \leq max \{C^{2}_{2}e^{-2t}, C^{2}_{3}e^{2 - 3\beta t} \}.$$
And since $2-3\beta < 0$. That implies that $(\Sigma T^{2}_{i}(t)e^{-2t}) \in L^{1}(\mathbf{R}^{+})$. 
 
 By assumption $\phi(t) \in L^{1}(\mathbf{R}^{+})$ that forces $(\Sigma T_{i}(t))' \in L^{1}(\mathbf{R}^{+})$. Which means $\Sigma T_{i}(t)$ is bounded. Again since each $T_{i}(t)$ is bounded from above (equation (*)), we conclude that the $|T_{i}(t)|$ 's are uniformly bounded for all $t$ large.

We have therefore established the inequalities for $S_{i}^{j}$.

\end{proof}
\begin{rem}
It would be interesting to know if one can still get these estimates, under weaker assumptions on the function $\phi(t)$, such as $\phi(t) \geq 0$ and $\lim_{t \to \infty} \phi(t) = 0$. It is quite clear that the Ricatti equation alone is not sufficient. To see this, if we take $\Sigma T_{i}(t) = \sqrt{t}$, we find  that (2.18) stays valid with $\phi$ decaying to zero at infinity.
\end{rem} 

\section{Scalar Curvature and Volume Bounds on Level Sets}
Next, we prove the following Lemma,

\begin{lem}\label{lem3.2.1}

 Let $(M^{n+1}, g)$, $n \geq 3$ be a complete conformally compact manifold with a pole $p$ and $Ric \geq -ng$. Let $t(x)= dist_{g}(p,x)$ and $\frac{1}{Vol(\Sigma_{t})}\int_{\Sigma_{t}}|R +n(n+1)| = o(e^{-2t})$, where for each $t$, $\Sigma_{t}$ is the level set of $t$, and $R$, the scalar curvature of $(M^{n+1}, g)$. The scalar curvatures and volumes of $(\bar{\Sigma}'_{t}, \bar{\gamma}'_{t})$ of $t$, satisfy the following inequalities,
 
 \begin{equation} 
 \int_{\bar{\Sigma}'_{t}}\bar{R}_{t}' \leq n(n-1)Vol(\bar{\Sigma}'_{t})+ o(1)
 \end{equation}
 
 \begin{equation} 
  Vol(\bar{\Sigma}'_{t}) \leq \omega_{n}
 \end{equation} 
where $\omega_{n}$ is the volume of $( S^{n}, \gamma_{0})$.  
 
 In particular, if $|R+n(n+1)| = o(e^{-2t})$ then one has $ \bar{R}'_{t} \leq n(n-1) + o(1)$.

\end{lem}

\begin{proof}

We will show that for the compactification $\bar{g}' = e^{-2t}\mathrm{g}$, 

$$ \int_{\bar{\Sigma}'_{t}}\bar{R}_{t}' \leq 4n(n-1)Vol(\bar{\Sigma}'_{t})+ o(1).$$
We refer to the following equations,

\begin{eqnarray}
\label{eq3.2.3}
\bar{Ric}' & = & Ric - (n-1)(Ddt - dt \circ dt) + (\Delta t - (n-1)|dt|^{2}) g, \\
\label{eq3.2.4}
\bar{R}' & = & e^{2t}(R + 2n \Delta t - n(n-1)|dt|^{2}).
\end{eqnarray}
From corollary ~(\ref{cor1}), we know that the second fundamental form $|\bar{S}'| \rightarrow 0 $ as $t \rightarrow \infty$.
By Gauss - Codazzi equation and corollary~(\ref{cor1}), we get:

\begin{equation}
\label{eq3.2.5}
  \bar{R}'_{t} = \bar{R}' - 2 \bar{Ric}'(\bar{N}', \bar{N}') + o(1).   
\end{equation}
Our assumption on the lower bound of the Ricci curvature gives the following upper bound on the Laplacian of $t$, $\Delta t \leq n\coth(t)$. Since $t$ is the distance function, we also have the following $g(\nabla t, \nabla t) = |dt|^{2} = 1 $
~(\ref{eq3.2.4}) becomes, 
\begin{equation}
\label{eq3.2.6}
  \bar{R}' = e^{2t}(R + 2n \Delta t - n(n-1)).
\end{equation}
Let $\bar{N}'$ be the unit normal to the level set $\bar{\Sigma}'_{t}$ in $\bar{g}'$ metric. Then from ~(\ref{eq3.2.3}) we get:  

\begin{eqnarray}
 \bar{Ric}'(\bar{N}', \bar{N}') & = & Ric(\bar{N}', \bar{N}') - (n-1)(Ddt - dt \circ dt)(\bar{N}', \bar{N}')  \nonumber \\
 & & +  (\Delta t - (n-1)) g(\bar{N}', \bar{N}') \nonumber \\
 & = & Ric(e^{t}N, e^{t}N) - (n-1)(Ddt - dt \circ dt)(e^{t}N, e^{t}N) \nonumber \\
 & & + (\Delta t - (n-1)) g(e^{t}N, e^{t}N) \nonumber \\
 & = & e^{2t}[Ric(N,N) - (n-1)(Ddt - dt \circ dt)(N,N) \nonumber \\
 &   & \quad + (\Delta t - (n-1))]  \nonumber
\end{eqnarray}   								
Since $Ddt(N,N) = g (\nabla^{2}_{N} t, N) = 0 $ and $dt \circ dt(N,N) = |dt|^{2} = 1 $, we get:

\begin{equation}
\label{eq3.2.7}
 \bar{Ric}'(\bar{N}', \bar{N}') = e^{2t} (Ric(N,N) + \Delta t )
\end{equation}
Substituting ~(\ref{eq3.2.6}) and ~(\ref{eq3.2.7}) in ~(\ref{eq3.2.5}), we get:  

  $$ \bar{R}'_{t} = e^{2t} [ R + 2(n-1) \Delta t - n(n-1) - 2Ric(N,N) ] + o(1). $$
Replacing $\Delta t $ by $n(\coth(t) - 1) + n$ on the right hand side, we get the following inequality:

\begin{eqnarray}
 \bar{R}'_{t} & \leq & e^{2t} [ R +2(n-1)(\frac{2ne^{-2t}}{1-e^{-2t}} + n) - n(n-1) - 2Ric(N,N) ]+ o(1) \nonumber  \\
  & \leq & e^{2t} [ R + \frac{4n(n-1)e^{-2t}}{1-e^{-2t}} + n(n-1) - 2Ric(N,N) ] + o(1) \nonumber  \\
  & \leq & e^{2t} [ R + n(n+1) + \frac{4n(n-1)e^{-2t}}{1-e^{-2t}} - 2(Ric(N,N) + n) ] + o(1) \nonumber
\end{eqnarray} 
Since $Ric(N,N) + n \geq 0$, we have,

 \begin{equation}
 \label{eq3.2.8}
 \bar{R}'_{t} \leq e^{2t}[ R + n(n+1) + \frac{4n(n-1)e^{-2t}}{1-e^{-2t}}] + o(1)
 \end{equation} 
Let $d \bar{\eta}'_{t}$ and $d\eta_{t}$ be the volume elements on $(\bar{\Sigma}'_{t}, \bar{g}'|_{\Sigma_{t}})$  and $(\Sigma_{t}, g|_{\Sigma_{t}})$ respectively. So, 
 
\begin{equation}
\label{eq3.2.10}
 d \bar{\eta}'_{t} = e^{-nt} d \eta_{t} 
\end{equation}  
Then,
 
\begin{eqnarray}
\label{eq3.2.11}
 \int_{\Sigma_{t}} \bar{R}'_{t} d \bar{\eta}'_{t} & \leq & e^{2t}\int_{\Sigma_{t}}[R + n(n+1) + \frac{n(n-1)e^{-2t}}{1-e^{-2t}}] d \bar{\eta}'_{t} + o(1)  \nonumber\\
  & \leq & e^{-(n-2)t }\Big ( \int_{\Sigma_{t}} |R+n(n+1)| d \eta_{t} \Big ) + \frac{n(n-1)}{1-e^{-2t}} Vol(\bar{\Sigma}'_{t}) + o(1) \nonumber \\
  &  &    
\end{eqnarray}  
By our assumption on the integral of the scalar curvature, we have:
 
 \begin{eqnarray}\label{eq3.2.12}
   e^{-(n-2)t}\Big ( \int_{\Sigma_{t}} |R+n(n+1)| d \eta_{t} \Big ) & = & e^{-(n-2)t} Vol(\Sigma_{t}). o(e^{-2t})  \nonumber\\
    & = & e^{-nt}Vol(\Sigma_{t}).o(1) \nonumber\\
    & = & Vol(\bar{\Sigma}'_{t}).o(1) \nonumber\\
    & = & o(1). 
 \end{eqnarray}
Substituting ~(\ref{eq3.2.12}) in ~(\ref{eq3.2.11}) we get:

\begin{equation}
\int_{\Sigma_{t}} \bar{R}'_{t} d \bar{\eta}'_{t} \leq \frac{n(n-1)}{1-e^{-2t}} Vol(\bar{\Sigma}'_{t}) + o(1). \nonumber
\end{equation}
Expanding the denominator in Taylor series, we see that:

\begin{equation}
\label{eq3.2.13}
\int_{\Sigma_{t}} \bar{R}'_{t} d \bar{\eta}'_{t} \leq n(n-1) Vol(\bar{\Sigma}'_{t}) + o(1).
\end{equation}
In particular, if $|R+n(n+1)| = o(e^{-2t})$, from ~(\ref{eq3.2.8}) we have:

$$ \bar{R}'_{t} \leq 4n(n-1) +o(1) $$
For the compactification $\bar{g}' = 4e^{-2t}g$, we have:

\begin{equation}
\bar{R}'_{t} \leq n(n-1) + o(1)
\end{equation}

\textbf{Proof of the volume bound}

From Gromov's relative volume comparison Theorem we have,

$$ t \rightarrow \frac{Vol(B(p, t)}{Vol(B_{0}(0, t))} $$
is a non-increasing function of $t$, where $B_{0}(0,t)$ is the ball of radius $t$ centered at the origin of the hyperbolic space, $(H^{n+1}, g_{-1})$. It also gives an upper bound on the ratio namely,

$$ \frac{Vol(B(p, t)}{Vol(B_{0}(0, t))}  \leq 1. $$
The volume elements of $(M,g)$ and $(H, g_{-1})$ in polar coordinates are  $\lambda(r,\theta)drd\theta$  and $\lambda_{-1}(r,\theta)drd\theta $, respectively. Therefore

\begin{eqnarray}
 Vol(B(p,t)) & = & \int_{0}^{t} \int_{S^{n-1}} \lambda(r,\theta)drd\theta,  \nonumber \\
 Vol(B_{0}(0,t)) & = & \int_{0}^{t} \int_{S^{n-1}} \lambda_{-1}(r,\theta)drd\theta.  \nonumber
\end{eqnarray}
So,

 $$ \frac{Vol(B(p, t)}{Vol(B_{0}(0, t))} = \frac{\int_{0}^{t} \int_{S^{n-1}} \lambda(r,\theta)drd\theta}{\int_{0}^{t} \int_{S^{n-1}} \lambda_{-1}(r,\theta)drd\theta}. $$
Since this is a non-increasing function, its derivative is non-positive, so the following inequality must hold

$$\frac{\int_{S^{n-1}} \lambda(t,\theta)d\theta}{\int_{S^{n-1}} \lambda_{-1}(t,\theta)d\theta} \leq \frac{\int_{0}^{t} \int_{S^{n-1}} \lambda(r,\theta)drd\theta}{\int_{0}^{t} \int_{S^{n-1}} \lambda_{-1}(r,\theta)drd\theta}. $$
Therefore,

$$ \frac{Vol(\Sigma_{t})}{Vol(S_{t})} \leq \frac{Vol(B(p, t)}{Vol(B_{0}(0, t))}  \leq 1, $$
where $S_{t}$ is the level set of the distance function from the origin at $t$ in the hyperbolic space.
So,

  $$ Vol(\Sigma_{t}) \leq Vol(S_{t}). $$
Or, equivalently
\begin{equation}
\label{eq3.2.14}
Vol(\bar{\Sigma}'_{t}) = \frac{Vol(\Sigma_{t})}{\frac{e^{nt}}{2^{n}}} \leq \frac{Vol(\Sigma_{t})}{\sinh^{n}t} \leq \frac{Vol(S_{t})}{\sinh^{n}t} = \omega_{n}
\end{equation}

\end{proof}

\section{Lipschitz Bound on Function \textit{c}}

On one hand, per our assumption $(\bar{M}, \bar{g})$ is a smooth compact manifold with $\bar{g}$ extending to a smooth metric up to the boundary $\partial M$. While on the other, $(\bar{M}, \bar{g}')$ is another compactification with $\bar{g}'$ smooth in $M$. We do not know yet, whether $\bar{g}'$ has a smooth extension on $\partial M$. In the next few Lemmas we show that in general the extension can only be $C^{0,1}$. But, if the conformal infinity is $(S^{n}, \gamma_{0})$ as in the Theorem, then the extension is $C^{\infty}$. As a first step we choose local coordinates in a neighborhood of the boundary as follows.

Since $\partial M $ is compact,it can be covered by finitely many coordinate n-balls $S_{k}(2R)$, such that $S_{k}(R)$ also cover $\partial M$. Each $S_{k}(2R)$ can be extended to $U_{k}(2R) = S_{k}(2R) \times [0, \epsilon_{0} ) \subset \bar{M}$, such that $\cup U_{k}(2R)$ and $\cup U_{k}(R)$ both cover a tubular neighborhood of $\partial M$ and there exist smooth coordinates in $U_{k}(2R)$, $\{ \theta_{1}, ... , \theta_{n}, \theta_{n+1} \}$ such that at any point $x \in M \cap U_{k}(2R)$,  $\bar{g}' = ds^{2} + \bar{g}'_{ij}d\theta^{i}d\theta^{j} $ where $s(x)=e^{-t(x)}$  and, $g = dt^{2} + e^{2t}\bar{g}'_{ij}d\theta^{i}d\theta^{j} = dt^{2} + g_{ij}d\theta^{i}d\theta^{j}$. Let us define: $S_{k;t}(2R) = U_{k}(2R) \cap \Sigma_{t}$ and $c_{t} = c|_{\Sigma_{t}}$, and $|\bar{\nabla}'^{t} c_{t}|^{2} = \bar{g}'^{ij}\partial_{i}c \partial_{j} c $

\begin{lem}\label{lem3.3.1}

  We have the following asymptotic bound on the $C^{0,1}$ norm of the function $c$: 

 (1)  $|\nabla c(x)|^{2} = O(e^{-2t})$

 (2)  $|\bar{\nabla}'c(x)|^{2} = e^{2t}|\nabla c(x)|^{2} = O(1)$
 
 Therefore,  $|\bar{\nabla}'^{t}c_{t}| < \Lambda_{1} $ for constant $\Lambda_{1}$ and on every  $S_{k;t}(2R)$ for $t$ large.
 
\end{lem}

\begin{proof}

We first observe that $g(\nabla c, \nabla c) = |\nabla r|^2 + |\nabla t|^2 - 2g(\nabla t, \nabla r) = 2(1 - g(\nabla r, \nabla t))$.
 
 Define $\phi_{1}(t) = e^{t}( 1 - g(\nabla r, \nabla t)) $. Clearly $\phi_{1}(t) \geq 0 $. 

 \begin{eqnarray}
 \partial_{t} \phi_{1}(t) & = & e^{t}(1 - g(\nabla r, \nabla t) ) - e^{t}\partial_{t} g(\nabla r, \nabla t) \nonumber\\
   & = & \phi_{1}(t) -  e^{t} [ g(\nabla_{\partial_{t}} \nabla r, \nabla t) + g(\nabla_{\partial_{t}} \nabla t, \nabla r)] \nonumber\\
   & = & \phi_{1}(t) -  e^{t} g(\nabla_{\partial_{t}} \nabla r, \nabla t) \nonumber\\
   & = & \phi_{1}(t) -  e^{t} \nabla^2 r(\nabla t, \nabla t). \nonumber
 \end{eqnarray}
 
  Let us write $ \nabla t = g(\nabla r, \nabla t) \nabla r + \sqrt{1 - g^{2}(\nabla r, \nabla t)} \nabla u $, where $\nabla u$ 
  is a unit vector orthogonal to $\nabla r$ . Making this substitution above, and using the asymptotic estimate of $|\nabla^2 r| $ from Lemma~(\ref{lem2}), we get,

 \begin{eqnarray}
\partial_{t} \phi_{1}(t) & = &  \phi_{1}(t) - e^{t} ( 1 - g^{2}(\nabla r, \nabla t)) \nabla^{2} r (\nabla u, \nabla u) \nonumber \\
  & = & \phi_{1}(t) - e^{t}(1 - g(\nabla r, \nabla t))(1 + g(\nabla r, \nabla t))(1 + O(e^{-\frac{3}{2}\beta t})) \nonumber \\
  & = & \phi_{1}(t) - \phi_{1}(t) (1 + g(\nabla r, \nabla t))(1 + O(e^{-\frac{3}{2}\beta t})). \nonumber
 \end{eqnarray} 
Since, $ 1 + g(\nabla r, \nabla t) \geq 1 $, we get the following inequality,

$$ \partial_{t} \phi_{1}(t) \leq \phi_{1}(t) - \phi_{1}(t)(1 + O(e^{-\frac{3}{2}\beta t})), $$
or,

$$ \partial_{t} \phi_{1}(t) \leq  \phi_{1}(t)(O(e^{-\frac{3}{2}\beta t})). $$
Integrating this, 
 \begin{equation} 
 \label{eq3.3.1}
   1 - g(\nabla r, \nabla t) = O(e^{-t}). 
\end{equation}
Next, we define another function $ \phi_{2}(t) = e^{2t}( 1 - g(\nabla r, \nabla t)) $.  
Following the same procedure, we get,
 \begin{equation}
 \label{eq3.3.2}
 \partial_{t} \phi_{2}(t) = 2\phi_{2}(t) - e^{2t}(1 - g(\nabla r, \nabla t))(1 + g(\nabla r, \nabla t))(1 + O(e^{-\frac{3}{2}\beta t})). 
 \end{equation}
As before, we replace $e^{2t}(1 - g(\nabla r, \nabla t))$ by $\phi_{2}(t)$ and use the above estimate, ~(\ref{eq3.3.2}) to get,
 
 $$ \partial_{t} \phi_{2}(t) = 2\phi_{2}(t) - \phi_{2}(t)(2 + O(e^{-t}))(1 + O(e^{-\frac{3}{2}\beta t})), $$
 or,
 
  $$ \partial_{t} \phi_{2}(t) = \phi_{2}(t)O(e^{-t}). $$
Therefore, we get
  
  $$ \frac{1}{2}|\nabla c|^{2} = 1 - g(\nabla r, \nabla t) = O(e^{-2t}) $$
That proves both (1) and (2). To prove the last claim we work with local coordinates,
 
  $$ |\bar{\nabla}'c (x)|^{2} = (\partial_{s} c(x))^{2} + \bar{g}^{ij}\partial_{i}c \partial_{j}(c) (x), $$
or,
  
  $$ |\bar{\nabla}'c (x)|^{2} = (\partial_{s} c(x))^{2} + |\bar{\nabla}'^{t}c_{t}(x)|^{2}. $$
Therefore, for some constant $\Lambda_{1}$
 
 \begin{equation}
 \label{eq3.3.3}
|\bar{\nabla}'^{t}c_{t}|^{2} \leq \Lambda_{1}.
\end{equation}
\end{proof} 
Under the additional assumption on the Ricci curvature in Lemma~(\ref{lem2}), we can prove the following Lemma. This result is not necessary for the proof of the main result.

\begin{lem}
There exists $\alpha \in (0,1) $ and a constant $\Lambda_{4}$ such that for all $i$ and $t$,
\begin{equation} 
 || c_{t} ||_{C^{1, \alpha}(S_{i;t}(R))} \leq \Lambda_{4}.
\end{equation}
\end{lem}

\begin{proof} 
We claim that, the scalar curvature $\bar{R}' $ is bounded from both above and below. This is true because, $\bar{R}' = e^{2t}(R + 2n\Delta t - n(n-1)) $ and we have already shown in Lemma (\ref{lem2}) that $ \Delta t = n + O(e^{-2t})$, and $ |R + n(n+1)| = o(e^{-2t})$. Therefore, we see that $|\bar{R}'|$ is bounded. Now we will show, $| \bar{\Delta}'c | = O(1) $.  We consider the two metrics $\bar{g}$ and $\bar{g} = e^{2c}\bar{g}'$. Therefore the scalar curvature of $(M,\bar{g})$ is given as,
 
  $$ \bar{R} = e^{-2c}( \bar{R}' + 2n \bar{\Delta}'c -n(n-1)|\bar{\nabla}'c|^{2})$$
Observe that all except $\bar{\Delta}'c$ are bounded, and hence $|\bar{\Delta}'c| = O(1)$. 
Next we show $ | \bar{\Delta}'^{t} c_{t} | = O(1)$.
In local coordinates, in each $U_{i}(2R) \cap M$, 
\begin{eqnarray}
  \lefteqn{ \bar{\Delta}'c = \frac{1}{\sqrt{det(\bar{g}')}} \big ( \partial_{i} \sqrt{det(\bar{g}')} \bar{g}'^{ij} \partial_{j} c \big ) + \frac{1}{\sqrt{det(\bar{g}')}} (\partial_{s} (\sqrt{det(\bar{g}')}\partial_{s} c))}, \nonumber \\
  &  &=  \bar{\Delta}'^{t} c_{t} + \partial_{s}(\partial_{s}c) + \frac{\partial_{s}\sqrt{det(\bar{g}')}}{ \sqrt{det(\bar{g}')}} \partial_{s} c. \nonumber
 \end{eqnarray} 
Therefore,
 
\begin{eqnarray}
  \lefteqn{ \partial_{s}(\partial_{s}c) = e^{2t} \big ((g(\nabla r, \nabla t) -1) + \partial_{t} g(\nabla r, \nabla t) \big )} \nonumber \\
  &  & = e^{2t} \big ((g(\nabla r, \nabla t) -1) +  \nabla^{2}r(\nabla t, \nabla t) \big ) = O(1). \nonumber
 \end{eqnarray} 
We showed in Lemma~(\ref{lem3.3.1}) that, $g(\nabla r, \nabla t) - 1 = O(e^{-2t})$ and by remark~(\ref{rem3.1.1}), $\nabla^{2}(\nabla t, \nabla t) = O(e^{-2r}) = O(e^{-2t})$. In addition, we have proved the following estimate, $\Delta t = n + O(e^{-2t})$ in Lemma~(\ref{lem2})we have, $\frac{\partial_{t}\sqrt{det(\bar{g}')}}{ \sqrt{det(\bar{g}')}} = \bar{\Delta}' s = e^{2t}(\Delta t - n) = O(1)$. Furthermore, $|\partial_{s} c| \leq |\bar{\nabla}' c|$ is uniformly bounded. Therefore, we get,
 $$ | \bar{\Delta}'^{t} c_{t} | = O(1). $$   
Finally, we apply Schauder estimates to show that for some constant $\Lambda_{4}$ and for all $i$, $t$ and $\alpha \in (0,1)$ 
 
 $$ || c_{t} ||_{C^{1, \alpha}(S_{i;t}(R))} \leq \Lambda_{3} \big ( ||\bar{\Delta}'^{t} c_{t} ||_{C^{0}(S_{i;t}(2R))} + ||c_{t}||_{C^{\alpha}(S_{i;t}(2R))} \big ) \leq  \Lambda_{4} $$ 
\end{proof}

\section{Regularity of the Boundary Metric}

 Thus we conclude that $ || c_{t} ||_{C^{1}(S_{k;t}(2R))} \leq \Lambda_{2} $, for some constant $\Lambda_{2}$ and for all $k$ and $t$. In other words, we found a Lipschitz bound on $c_{t}$. Hence, by Arzela-Ascoli's Theorem, we conclude that given a sequence $( \bar{\Sigma}_{t_{i}}', \bar{\gamma_{t_{i}}}' = \bar{g}'|_{\bar{\Sigma}_{t_{i}}'} )$, there is a subsequence that converges to $(\Sigma, \gamma)$, where $\gamma$ is $C^{0,1}\cap W^{1,p}$ conformal to $\gamma_{0}$, for any $p \geq 1$.

So far we have not used the structure of the boundary explicitly, that is all our results hold under the assumption that  boundary $(\partial M, [\bar{g}|_{\partial M}])$ is smooth, not necessarily $(S^{n}, [\gamma_{0}])$ . Now for the first time we will use the assumption that the boundary is $(S^{n}, [\gamma_{0}])$ to show that $\gamma$ is in fact smooth. Now that we have $\gamma$, a $C^{0,1}\cap W^{1,p}$ metric, for any $p \geq 1$, we can define scalar curvature of $(\Sigma, \gamma)$, in the sense of distributions. We choose a representative element from the conformal class $[\gamma_{0}]$. Without any loss of generality we pick the round metric $\gamma_{0}$. To simplify our computations below, we define $\gamma = u^{\frac{4}{n-2}} \gamma_{0} $, where $u$ is $C^{0,1}\cap W^{1,p}$,for all $p \geq 1$.

\begin{lem}
If $(\partial M, [\gamma]) = (S^{n}, [\gamma_{0}])$, where $\gamma = u^{\frac{4}{n-2}} \gamma_{0}$ for some $u \in W^{1,p}(S^{n})$, for any $p \geq 1$, then $\gamma$ is smooth. 
\end{lem}
\begin{proof}
 The scalar curvature of $(\Sigma, \gamma)$, in the weak sense is given by,
 
 \begin{equation}
 \label{eq3.4.1}
 \bar{R}' = u^{-\frac{n+2}{n-2}} \big ( n(n-1)u - \frac{4(n-1)}{n-2}\Delta_{\gamma_{0}} u \big ),
 \end{equation}
where $\Delta_{\gamma_{0}}u$ is in $W^{-1,p}$, since $\Delta_{\gamma_{0}}: W^{1,p} \rightarrow W^{-1,p}$. Therefore we can integrate $\bar{R}'$.
Let $d\eta$ and $d\eta_{0}$ be the volume elements of $(\Sigma, \gamma)$ and $(\Sigma, \gamma_{0})$ respectively. Therefore,
  $d\eta = u^{\frac{2n}{n-2}}d\eta_{0}$.
Integrating the scalar curvature over $\Sigma = S^{n}$,
 
 \begin{equation}  
 \label{eq3.4.2}
 \int_{S^{n}} \bar{R}' u^{\frac{2n}{n-2}}d\eta_{0} = \int_{S^{n}} \big( n(n-1)u^{2} +  \frac{4(n-1)}{n-2}   |\nabla_{\gamma_{0}}u|^{2} \big ) d\eta_{0}
 \end{equation}
Recall, we showed earlier in Lemma~(\ref{lem3.2.1}) that the scalar curvatures of level sets of $t$, 
\begin{eqnarray}
 \bar{R}_{t} & \leq & n(n-1) + o(1), \nonumber \\
 Vol(\bar{\Sigma}_{t}') & \leq & \omega_{n}. \nonumber 
\end{eqnarray} 
 Therefore, at the boundary we have 
 \begin{eqnarray}
 \int_{S^{n}} \bar{R}' u^{\frac{2n}{n-2}}d\eta_{0} & \leq & n(n-1) \int_{S^{n}}u^{\frac{2n}{n-2}}d\eta_{0}, \nonumber \\
 \int_{S^{n}}u^{\frac{2n}{n-2}}d\eta_{0}& \leq & \omega_{n}. \nonumber 
\end{eqnarray}
So,
  
 \begin{eqnarray}  
 n(n-1)\int_{S^{n}} u^{\frac{2n}{n-2}}d\eta_{0} & \geq &  \int_{S^{n}} \big ( n(n-1)u^{2} +  \frac{4(n-1)}{n-2}|\nabla_{\gamma_{0}}u|^{2} \big ) d\eta_{0}, \nonumber \\
n(n-1)\big ( \int_{S^{n}} u^{\frac{2n}{n-2}}d\eta_{0} \big )^{\frac{2}{n}} & \geq &  \frac{\int_{S^{n}} \big ( n(n-1)u^{2} +  \frac{4(n-1)}{n-2} |\nabla_{\gamma_{0}}u|^{2} \big ) d\eta_{0}}{\big ( \int_{S^{n}} u^{\frac{2n}{n-2}}d\eta_{0} \big )^{\frac{n-2}{n}}}.   \nonumber
 \end{eqnarray} 
Replacing  $\int_{S^{n}}u^{\frac{2n}{n-2}}d\eta_{0}$ by $\omega_{n}$, we get

\begin{equation}
\label{eq3.4.3}
 n(n-1)(\omega_{n})^{\frac{2}{n}} \geq  \frac{\int_{S^{n}} \big ( n(n-1)u^{2} +  \frac{4(n-1)}{n-2} |\nabla_{\gamma_{0}}u|^{2} \big ) d\eta_{0}}{\big ( \int_{S^{n}} u^{\frac{2n}{n-2}}d\eta_{0} \big )^{\frac{n-2}{n}}}.
\end{equation} 
 
The right-hand side is the Yamabe quotient, $Q(\gamma)$ on $S^{n}$. So, $\mathrm{inf} \{ Q(\gamma) | \gamma \in [\gamma_{0}] \} = n(n-1) (\omega_{n})^{\frac{2}{n}} $. This can be easily shown by taking the stereographic projection $\pi: S^{n} \rightarrow \mathbf{R}^{n+1}$ and writing out the right-hand side in $\mathrm{R}^{n}$ with respect to the pulled-back metric, $\pi^{*} \gamma_{0}$, and comparing the integral with the Sobolev constant.ting out the right-hand side in $\mathrm{R}^{n}$ with respect to the pulled-back metric, $\pi^{*} \gamma_{0}$, and comparing the integral with the Sobolev constant.

Therefore, we conclude that the above inequality is an equality, that is $\bar{R}' = n(n-1)$ (a.e.) on $S^{n}$. But that would imply that $\bar{R}' \in W^{1,p}(S^{n})$ for all $p > 1$. This implies, by Sobolev embedding, that $\bar{R}'$ is continuous on $S^{n}$. Clearly, the subset $\mathcal{U}$ of the boundary $S^{n}$ where $\bar{R}' = n(n-1)$ is a set of full measure. Therefore for any $x \in S^{n} - \mathcal{U}$, we can find a sequence of points in $\mathcal{U}$ converging to $x$. Thus, by continuity of $\bar{R}'$, $x$ is in $\mathcal{U}$, and therefore $\mathcal{U} = S^{n}$. Thus,

\begin{equation}
 \bar{R}' = n(n-1) 
\end{equation} 
 
Finally, a standard boot-strapping argument applied to the semi-linear elliptic P.D.E., 

 \begin{equation}
 \bar{R}'= n(n-1) = u^{-\frac{n+2}{n-2}} \big ( n(n-1)u - \frac{4(n-1)}{n-2}\Delta_{\gamma_{0}} u \big ) 
 \end{equation}
shows that $u$ is in fact smooth. 
\end{proof} 
 
 Hence we have shown that the two metrics $\gamma$ and $\gamma_{0}$ on $\Sigma$ are smoothly conformal to one another and hence $\gamma$ is smoothly conformal to the round metric. Applying Obata's Theorem \cite{ob} we conclude that $( \Sigma, \gamma)$ is isometric to $(S^{n}, \gamma_{0})$.   
 
\section{Proof of the Main Theorem}
 
 \begin{proof}
 
 We proved the following inequality earlier, ~(\ref{eq3.2.14})

  $$ Vol(\bar{\Sigma}'_{t}) = \frac{Vol(\Sigma_{t})}{\frac{e^{nt}}{2^{n}}} \leq \frac{Vol(\Sigma_{t})}{\sinh^{n}t} \leq \frac{Vol(S_{t})}{\sinh^{n}t} = \omega_{n}. $$
Furthermore, we have just shown, that the level sets $(\bar{\Sigma}'_{t}, \bar{\gamma}'_{t})$ converge to $(S^{n}, \gamma_{0})$ at the boundary, therefore $ Vol(\bar{\Sigma}'_{t}) \rightarrow \omega_{n}$ as $ t \rightarrow \infty $.
So, combining the two and the fact that volume ratio is strictly non increasing, we see that,

\begin{center}
$ \frac{Vol(B(p, t)}{Vol(B_{0}(0, t))}  \equiv 1 $   for all $t$.
\end{center}
This implies that $ \Delta t = n\coth (t) $.
Therefore,
 
 \begin{equation}
 -n \geq \partial_{t}(\Delta t) + | \nabla^{2}t|^{2} \geq \partial_{t} (\Delta t) + \frac{ (\Delta t)^{2}}{n} = -n. 
 \end{equation}
The first inequality comes from the Ricatti equation and $Ric \geq -ng $.
This implies that:

\begin{equation}
   (\Delta t)^{2} = n|\nabla^{2}t|^{2}.
\end{equation}
So, if $\{ \lambda_{1}, ... \lambda_{n} \}$ are the eigenvalues of $\Delta^{2} t$ then the above equality implies: 
\begin{equation}
(\lambda_{1} + ... + \lambda_{n})^{2} = n(\lambda_1^{2} + ... + \lambda_{n}^{2}).
\end{equation}
This can happen if and only if $\lambda_{1} = \lambda_{2} = ... = \lambda_{n} = \lambda $. Which means that 
 $$ \Delta t = n \lambda = n \coth(t), $$
or,
 $$ \lambda = \coth(t). $$
This, in turn, implies that the metric $g$ on $M$ is a space form and is precisely the hyperbolic metric, $g = dt^{2} + \sinh^{2}(t) \gamma_{0} $.     
That completes the proof of the main Theorem. 

\end{proof}

\begin{flushleft}
Department of Mathematics \\
SUNY at StonyBrook, NY 11794 \\
email: sunny@math.sunysb.edu
\end{flushleft}

\end{document}